\def\bu{\bullet}
\def\marker{\>\hbox{${\vcenter{\vbox{
    \hrule height 0.4pt\hbox{\vrule width 0.4pt height 6pt
    \kern6pt\vrule width 0.4pt}\hrule height 0.4pt}}}$}\>}
\def\gpic#1{#1
     \medskip\par\noindent{\centerline{\box\graph}} \medskip}
\documentclass[12pt]{article}

\usepackage[margin=2cm]{geometry}

\usepackage{amsmath, amssymb, latexsym, amsthm, setspace}

\newtheorem{theorem}{Theorem} 

\newtheorem{lemma}[theorem]{Lemma}

\newtheorem*{theorem*}{Theorem}
\newtheorem*{conjecture*}{Conjecture}
\newtheorem*{corollary*}{Corollary}
\newtheorem*{proposition*}{Proposition}

\theoremstyle{definition}

\theoremstyle{remark}

\renewcommand{\labelenumi}{(\alph{enumi})} 

\def\st{\colon\,}
\def\VEC#1#2#3{#1_{#2},\ldots,#1_{#3}}
\def\bT{{\bf T}}
\def\la{\langle}
\def\ra{\rangle}
\def\qed{\hfill$\Box$}

\usepackage[dvips]{graphicx}
\newcommand{\figref}[1]{Figure~\ref{#1}}
\graphicspath{{counterfigs/}}
\newlength{\mydim}

\author{
Armen S. Asratian\thanks{
Link\"oping University, Link\"oping Sweden, arasr@mai.liu.se.}\,,
Carl Johan Casselgren\thanks{
Ume\aa\enskip University, Ume\aa\,, Sweden, carl-johan.casselgren@math.umu.se.}
}

\title{On path factors  of $(3,4)$-biregular
bigraphs}
\date{}

\begin{document}
\maketitle

\bigskip 
{\bf Abstract}.  A $(3,4)$-biregular bigraph $G$ is
 a bipartite graph where all vertices in one part have  degree 3 and
all vertices in the 
other part 
have  degree 4.  A path factor of  $G$ 
 is a spanning subgraph whose components are  nontrivial
paths. We prove that 
a simple $(3,4)$-biregular bigraph 
always has a path 
factor such that the endpoints of each path have degree  three.
Moreover we suggest a 
polynomial
algorithm for the construction of such a path factor.\bigskip


\section{ Introduction}

We use \cite{graphtheory} and \cite{bipartite} for terminology and
notation not defined 
here and
consider  finite loop-free graphs only. $V(G)$ and $E(G)$ denote the
sets
of vertices and 
edges of a graph $G$, 
respectively. A \textit{proper edge coloring} of a graph $G$ with
colors $1,2,3,\dots$ is
a mapping $f: E(G)\rightarrow\{1, 2, 3, \dots\}$ such that 
$f(e_1)\not=f(e_2)$ for every pair  
of
adjacent edges $e_1$ and $e_2$.
A  bipartite graph with bipartition $(Y,X)$ is called
an \textit{$(a,b)$-biregular} bigraph 
if every vertex in $Y$ has degree $a$ and every vertex in $X$ has
degree $b$.
A {\it path factor} of  a graph $G$ 
 is a spanning subgraph whose components are nontrivial
paths. Some results on different types of path factors can be found 
in 
\cite{Aki,Ando, Kaneko, cubic2,Plummer,wang}.
In particular, Ando et al \cite{Ando} showed that a claw-free graph
with minimum degree 
$d$ has a path factor whose components are paths of length at least
$d$. Kaneko \cite
{Kaneko} showed that every cubic graph has a path factor such that
each component is a 
path of length $2,3$ or $4$. It was shown in \cite{cubic2} that  a
2-connected cubic graph  
has a path factor whose components are paths of length 2 or 3.  

 In this paper we investigate the existence of path factors of
$(3,4)$-biregular bigraphs such that the endpoints of each path have
degree  three. Our 
investigation is motivated by a problem on interval colorings.  A
proper edge coloring  of a graph $G$ with colors $1, 2, 3,\dots$ is
called
an {\it interval} (or {\it consecutive}) coloring if the
colors received by the edges
incident with each vertex  of $G$ form an interval of integers.
The notion of interval colorings was introduced in 1987 by Asratian
and
Kamalian  \cite{Kamalian1} (available in English as  \cite{Kamalian}).
Generally, it is an $\mathcal{NP}$-complete problem to determine
whether a given 
bipartite
graph has an interval coloring \cite{sevastjanov}.  Nevertheless,
  trees, 
regular and complete
bigraphs  \cite{hansen,Kamalian2}, doubly 
convex bigraphs \cite{Kamalian2}, grids \cite{timetabling4} and all 
outerplanar 
bigraphs \cite{Axe, timetabling3}
 have interval 
colorings.
Hansen  \cite{hansen} proved  that every
$(2,\beta)$-biregular bigraph admits an interval coloring if $\beta$
is an even
integer. A similar result for  $(2,\beta)$-biregular bigraphs for odd
$\beta$ was given in  \cite{2-biregular, unpublished}.
Only a little is known about $(3,\beta)$-biregular bigraphs. It
follows from the result of 
Hanson and Loten  \cite{LowerBound} that no 
 such a  graph   has an interval coloring with fewer than 
$3+b-\gcd(3,b)$ colors, where $
\gcd$ denotes the greatest common divisor.  We showed in \cite{biregular} that the
problem to determine whether a
$(3,\beta)$-biregular bigraph has an interval
coloring is $\mathcal{NP}$-complete in the case when 3
divides $\beta$. 

 It is unknown whether all $(3,4)$-biregular bigraphs have interval
colorings.
 Pyatkin \cite{CubicSub} showed that such a graph $G$  has an interval
coloring if $G$ 
has a 3-regular subgraph covering the vertices of degree four. 
Another sufficient condition for the existence of an interval
coloring of a $(3,4)$-biregular bigraph $G$ was obtained in
\cite{West,master}: $G$ admits an 
interval coloring if it has a path factor  where every component  is
a path of length not
exceeding 8 and the endpoints of each path have degree  three.  It was conjectured in 
\cite{West} that every simple $(3,4)$-biregular bigraph
has  such a path factor. However  this seems  difficult to prove. 

In this note we prove a little weaker result.  We show that 
a simple $(3,4)$-biregular bigraph 
always has a path 
factor  such that the endpoints of each path have degree three.
Moreover, we suggest a 
polynomial
algorithm for the construction of such a path factor.

Note that  $(3, 4)$-biregular bigraphs with multiple edges need 
not have  path factors with  the  required property. For example,
consider the 
graph $G$ formed from three triple-edges by adding a claw; that is,
the pairs $x_iy_i$ 
have multiplicity 
three for $i\in \{1, 2, 3\}$, and there is an additional vertex $y_0$
with neighborhood $\{x_1 , x_2 , x_3 \}$. 
Clearly, there is no path factor of $G$ such that the endpoints of
each path have degree 
3.

\section{The result}
\label{sec:factors}

A {\it pseudo path factor} of a $(3,4)$-biregular bigraph
$G$ with bipartition $(Y,X)$ is a  subgraph $F$ of $G$, such that
every 
component of $F$ is a  path  of  even length and $d_{F}(x) = 2$ for
every $x \in X$.  Let 
$V_F=\{y\in Y :  d_F(y)>0\}$. 

\begin{theorem}
\label{th:pseudofactor}
Every simple $(3,4)$-biregular bigraph  has a pseudo path factor.
\end{theorem}
\begin{proof}[\emph{\textbf{Proof.}}]
Let $G$ be a simple $(3,4)$-biregular bigraph with bipartition
$(Y,X)$.
The algorithm below   constructs  a sequence of
subgraphs $F_0,$$ F_1,$$ F_2, \dots$ of $G$,
where $V(F_0)=V(G)$,  $\emptyset=E(F_0)\subset E(F_1)\subset 
E(F_2)\subset \dots$ and each component of $F_j$ is a path, for every
$j\geq 0$.
At each step $i\geq 1$ the algorithm constructs $F_i$ by adding
to $F_{i-1}$ one or two edges until the condition $d_{F_j}(x)=2$ 
holds for all $x\in
X$,  where $j\geq 1$.  
Then $F=F_j$ is a pseudo path factor of $G$. Parallelly the algorithm 
constructs a
sequence of subgraphs 
$U_0, U_1, U_2, \dots$ of $G$, where $V(U_0)=V(G)$,
$\emptyset=E(U_0)\subset
E(U_1)\subset  
E(U_2)\subset \dots \subset E(U_j)$. The edges of each $U_i$ will not 
be in the final 
pseudo path factor $F$. 
The algorithm is based on Properties 1-4. During the algorithm the 
vertices in the set $Y$ 
are considered to
be unscanned or scanned. Initially all vertices in $Y$ are unscanned.
At the beginning of each step $i\geq 1$ we have a current vertex
$x_i$. The algorithm  
selects an unscanned vertex $y_i$, adjacent to $x_i$, and determines
which edges 
incident with $y_i$ will be in $F_i$ and which ones 
in $U_i$. If $d_{F_i}(v)=2$ for each $v\in X$, the algorithm stops.
Otherwise   the 
algorithm selects a new current vertex and goes to  the next step.
\\
\\
{\bf Algorithm} \medskip

Initially $F_0=(V(G),\emptyset)$,  $U_0=(V(G), \emptyset)$ and   all
vertices in $Y$ are 
unscanned. \medskip

\noindent	{\bf Step $0$}. Select a vertex $y_0\in Y$. Let $x_0, x_1,
w$
be the vertices in $X$ adjacent to $y_0$ in $G$. Put $F_1=F_0+\{wy_0,
y_0x_0\}$  and 
$U_1=U_0+y_0x_1$. 
Consider the vertex $y_0$  to be scanned. Go to step $1$ and consider
the vertex
$x_1$ as the current vertex for  step 1.\bigskip

\noindent {\bf Step $i$ $(i \geq 1)$.}  Suppose that a vertex $x_i$
with
$d_{F_{i-1}}(x_i)\leq 1$ was selected
at step $(i-1)$ as
the current  vertex. By Property 4 (see below), $d_{U_{i-1}}(x_i)\leq
2$. Therefore  there is an edge $x_iy_i$ with $y_i \in Y$ which
neither belongs to $F_{i-1}
$,
nor to $U_{i-1}$. Then, by Property 3, the vertex $y_i$  is an
unscanned vertex and 
therefore 
the subgraph $F_{i-1}
+x_iy_i$ does not contain a cycle.  Since $d_G(y_i)=3$, the vertex
$y_i$, besides $x_i$, is adjacent to two other
vertices, $w_1^{(i)}$ and $w_2^{(i)}$. \smallskip
	
\noindent \textbf{\textit{Case 1}}.
$d_{F_{i-1}}(w_1^{(i)})=2=d_{F_{i-1}}(w_2^{(i)})$.
	\\
	 Put $F_i=F_{i-1}+x_iy_i$ and 
$U_i=U_{i-1} + \{y_iw_1^{(i)}, y_iw_2^{(i)}\}$. Consider the vertex
$y_i$ to be scanned.
If $d_{F_i}(v)=2$ for every vertex $v\in X$ then Stop.  Otherwise
select an arbitrary vertex $x_{i+1}\in X$
with $d_{F_i}(x_{i+1})\leq 1$, go to step $(i+1)$ and consider
$x_{i+1}$ as
the current vertex for step $(i+1)$. \smallskip

\noindent \textbf{\textit{Case 2}}. $d_{F_{i-1}}(w_1^{(i)})=2$ and
$d_{F_{i-1}}(w_2^{(i)})\leq 1$. 
\\
Put $F_i=F_{i-1}+x_iy_i$, 
$U_i=U_{i-1}+\{y_iw_1^{(i)}, y_iw_2^{(i)}\}$ and consider the vertex $y_i$ to be scanned.
Furthermore put  $x_{i+1}=w_2^{(i)}$, go to step $(i+1)$ and consider
the vertex $x_{i+1}$ as the current  vertex  for step $(i+1)$. \smallskip

\noindent \textbf{\textit{Case 3}}.  $d_{F_{i-1}}(w_1^{(i)})\leq 1$ and
$d_{F_{i-1}}(w_2^{(i)})\leq 1$. 

\noindent \textsl{Subcase 3a}. $d_{F_{i-1}}(w_1^{(i)})=0$ or
$d_{F_{i-1}}(w_2^{(i)})=0$. 
\\
We assume that 
$d_{F_{i-1}}(w_1^{(i)})=0$. Put  $F_i=F_{i-1}+\{x_iy_i,
y_iw_1^{(i)}\}$,
$U_i=U_{i-1}+y_iw_2^{(i)}$ and consider the vertex $y_i$ to be scanned. 
Furthermore put   $x_{i+1}=w_2^{(i)}$, go to step $(i+1)$ and
consider the vertex $x_{i+1}$ as the current  vertex  for step $(i+1)$.

\noindent \textsl{Subcase 3b}. 
$d_{F_{i-1}}(w_1^{(i)})=1=d_{F_{i-1}}(w_2^{(i)})$. 
\\
Since $y_i$ is an unscanned vertex and $F_{i-1}+x_iy_i$ does not
contain a cycle, the 
vertex $y_i$  is an endvertex of only one path in $F_{i-1}+x_iy_i$. 
Then at least one of the graphs
$F_{i-1}+\{x_iy_i, y_iw_1^{(i)}\}$ and $F_{i-1}+\{x_iy_i,
y_iw_2^{(i)}\}$ does not contain a 
cycle. Assume, for example,  that $F_{i-1}+\{x_iy_i, y_iw_1^{(i)}\}$
does not contain a cycle. Then put $F_i=F_{i-1}+\{x_iy_i,y_iw_1^{(i)}\}$,
$U_i=U_{i-1}+y_iw_2^{(i)}$ and consider the vertex $y_i$ to be scanned.
 Furthemore put  $x_{i+1}=w_2^{(i)}$, go to step $(i+1)$ and
consider the vertex 
$x_{i+1}$ as the current  vertex  for step $(i+1)$.\bigskip

Now we will prove the correctness  of the algorithm. At the beginning
of  step $i$  we have that 
$x_i$ is the  current vertex, $y_i$ is an
unscanned vertex adjacent to $x_i$ and $w_1^{(i)}$, $w_2^{(i)}$ are
the two other vertices adjacent to $y_i$. The following  two 
properties are evident. \medskip

\noindent {\bf Property 1}.
The algorithm  determines which edges incident with
$y_i$ will be in 
$F_i$ and which edges will  be in $U_i$.  The vertex $y_i$ is then considered to be 
scanned and  the algorithm  will never  consider $y_i$ again. \medskip

\noindent {\bf Property 2}. The current vertex $x_{i+1}$ for step $(i+1)$ is
selected among the
vertices $ w_1^{(i)}$ 
and $w_2^{(i)}$,
except the case $d_{F_i}(w_1^{(i)})=d_{F_i}( w_2^{(i)})=2$ when an
arbitrary vertex $x_{i
+1}\in X$
with $d_{F_i}(x_{i+1})\leq 1$ is selected as the current vertex.
\bigskip

Properties 1 and 2 imply  the next property: \medskip

\noindent {\bf Property 3}.  If $x\in X$, $y\in Y$ and the edge $xy$ neither
 belongs to $F_{i-1}$, nor  to $U_{i-1}$, 
then
the vertex $y$ is unscanned at the beginning of  step $i$.
\\ \\
{\bf Property 4}. If $x\in X$ and $d_{F_{i-1}}(x)\leq 1$ then
$d_{U_{i-1}}(x)\leq 2$.

\begin{proof}[\emph{\textbf{Proof.}}]
The statement is evident if $d_{U_{i-1}}(x)=0$. Suppose
that $d_{U_{i-1}}(x)\geq 1$
and  $j$ is the minimum number such that $j<i$ and an  edge  incident
with $x$  was included in  $U_{j-1}$ at step 
$(j-1)$. Then the statement of Property 4 is evident if $j=i-1$.

Now we consider the case $j<i-1$.  Clearly, $d_{F_{j-1}}(x)\leq 1$
 because $F_{j-1}\subset F_{i-1}$ and $d_{F_{j-1}}(x)\leq 
d_{F_{i-1}}(x)\leq 1$.
Let $xy_{j-1}$ be the edge  included in  $U_{j-1}$ at step 
$(j-1)$.  Since $d_{U_{j-1}}(x)=1$
and $d_{F_{j-1}}(x)\leq 1$,  there is an edge $xy_j$ with $y_j\in Y$
which neither
belongs to $F_{j-1}$ , nor to
$U_{j-1}$. Then, by Property 3, the vertex $y_j$  is an unscanned
vertex and therefore  
the subgraph $F_{j-1}+xy_j$ does not contain a cycle. According to the
description of the 
algorithm, the edge $xy_j$ will be in any case included in 
$F_{j}$ at step $j$, that is, $d_{F_j}(x)\geq 1$.  Then $d_{F_k}(x)=1$
for every $k$, $j\leq k\leq i-1$, because $F_j\subset F_k\subset
F_{i-1}$ and $1\leq d_
{F_j}(x)\leq
d_{F_k}(x) \leq d_{F_{i-1}}(x)\leq 1$.
Now we will show that $d_{U_{k-1}}(x)=1$ for each $k$, $j\leq k<i-1$.
Suppose to the contrary  that $d_{U_{k-2}}(x)=1$ and 
$d_{U_{k-1}}(x)=2$ for some $k$,
$j<k<i-1$, that is, another edge  incident with $x$  was included in
$U_{k-1}$ at
step $(k-1)$. Then the conditions $d_{U_{k-1}}(x)=2$ and
$d_{F_{k-1}}(x)=1$ imply that 
there is an edge $e\not=y_jx$ 
incident with $x$ which neither belongs to $F_{k-1}$ , nor to
$U_{k-1}$.
Using a similar argument as above we obtain that the edge $e$ should
be included in 
$F_k$ at step $k$.
But then $d_{F_{i-1}}(x)\geq d_{F_{k}}(x)=2,$ which
contradicts our assumption $d_{F_{i-1}}(x)\leq 1$. Thus
$d_{U_{k-1}}(x)=1$
for each $k$, $j\leq k<i-1$.  It is possible that an edge incident
with $x$ will be included
in $U_{i-1}$ at step $(i-1)$. Therefore $d_{U_{i-1}}(x)\leq 2$.
\end{proof}

The description of the algorithm and Properties 1-4 show that the
algorithm will stop at step $i$ only when $d_{F_{i}}(x)=2$ for every
$x\in X$, that is,
when $F_{i}$ is a pseudo path factor of $G$. The proof of Theorem
\ref{th:pseudofactor}
is complete.
\end{proof}

Now we will prove that every pseudo path factor of a $(3,4)$-biregular
bigraph $G$ can be transformed into a path factor of $G$, such that the endpoints of each path
have degree 3.

\begin{lemma}
\label{lem:numbervertices}
	Let $G$ be a $(3,4)$-biregular bigraph with bipartition $(Y,X)$.
	Then
	$|X| = 3k$ and $|Y| = 4k$, for some positive integer $k$.
\end{lemma}

This is evident because $|E(G)|=4|X|=3|Y|$. 

\begin{lemma}
\label{lemmapseudo}
Let $F$ be a  pseudo path factor of a $(3,4)$-biregular bigraph $G$
with bipartition $(Y,X)$. Then $F$ has a component which is a path of
length at least four.
\end{lemma}
\begin{proof}[\emph{\textbf{Proof.}}]
By Lemma \ref{lem:numbervertices} we have that $|X|=3k$ and $|Y|=4k$
for some integer
$k$. We also have that $d_F(x)=2$ for each vertex $x\in X$. If the
length of all paths in $F
$
is two, then $|Y|\geq 2|X|=6k$ which contradicts $|Y|=4k$. Therefore 
$F$ has a 
component which is
a path of length at least four.
\end{proof}

\begin{theorem}
\label{th:pseudo}
Let $F$ be a  pseudo path factor of a simple $(3,4)$-biregular bigraph
$G$ with
bipartition $(Y,X)$. If $V_F\not=Y$ and $y_0$ is a vertex  with 
$d_F(y_0)=0$, then there is a pseudo path factor $F'$ with
$V_{F'}=V_F\cup \{y_0\}$, such that no path in $F'$ is longer than the longest path in $F$.
\end{theorem}

\begin{proof}[\emph{\textbf{Proof.}}]
Let $y_0\in Y$ and $d_F(y_0)=0$. We will describe an algorithm which
will 
construct a special trail $T$ with origin $y_0$.  \smallskip

\noindent {\bf Step $1$}.  Select an edge $y_0x_1\notin E(F)$. Since
$d_F(x_1)=2$, 
there are two edges of $F$, $x_1y_1$ and $x_1u_1$, which are incident 
with $x_1$. 

\noindent \textbf{\textit{Case 1}}.  $d_F(y_1)=2$ or $d_F(u_1)=2$.
\\
Suppose, for example,  that $d_F(y_1)=2$. Then put
$T=y_0\rightarrow x_1\rightarrow y_1$
and Stop. 

\noindent \textbf{\textit{Case 2}}.  $d_F(y_1)=1=d_F(u_1)$.
\\
Put $T=y_0\rightarrow x_1\rightarrow y_1$ and 
go to Step 2. \medskip

\noindent {\bf Step $i$ $(i \geq 1)$.}  Suppose that we have already
constructed a
trail $T=y_0\rightarrow x_1\rightarrow y_1\rightarrow \dots
\rightarrow
x_i\rightarrow y_i$ which satisfies the following conditions:

(a)  All edges in $T$ are distinct and $y_{j-1}x_{j}\notin E(F)$, 
$x_jy_j\in E(F)$   for $j=1,\dots ,i$.

(b) The vertices $y_1,\dots ,y_i$ are distinct.

(c) A component of $F$ containing the vertex $x_j$ is a path of length
2, for $j=1, \dots ,i$.

Select an edge $e\in E(G)\setminus E(F)$ which is incident with $y_i$.
The
existence of such an
edge follows from the conditions (a), (b) and (c). Moreover, the
condition (b) implies that
$e\notin T$. Let $e=y_ix_{i+1}$. Then $d_F(x_{i+1})=2$ because $F$ is
a  pseudo path 
factor of $G$.
Since $e\notin E(T)$, the conditions (a), (b) and (c) imply that at
least one of the edges
of $F$ incident with $x_{i+1}$, does not belong to $T$. 

\noindent \textbf{\textit{Case 1}}. $x_{i+1}$ lies on a component of $F$ which
is a path of
length two. 
\\ 
Select a vertex $y_{i+1}$ such that $x_{i+1}y_{i+1}\in E(F)\setminus
E(T)$,
add the edge  $x_{i+1}y_{i+1}$ and the vertex $y_{i+1}$ to $T$  and go
to step $(i+1)$. 
Now
$T=y_0\rightarrow x_1\rightarrow y_1\rightarrow \dots \rightarrow
x_{i+1}\rightarrow y_{i+1}.$

\noindent \textbf{\textit{Case 2}}. $x_{i+1}$ lies on a component of $F$ which
is a path of
length at least four. 
\\
There is a vertex $y_{i+1}$ such that $x_{i+1}y_{i+1}\in E(F)\setminus
E(T)$
and $d_F(y_{i+1})=2$. Add the edge  $x_{i+1}y_{i+1}$
and the vertex $y_{i+1}$ to $T$ and Stop. We have now that
$T=y_0\rightarrow x_1
\rightarrow 
y_1\rightarrow \dots \rightarrow
x_{i+1}\rightarrow y_{i+1}$. \medskip

By Lemma \ref{lemmapseudo}, $F$ has a component which is a path of
length at least
four. Therefore the 
algorithm will stop after a finite number of steps. Let the trail
$T=y_0\rightarrow x_1\rightarrow
y_1\rightarrow \dots \rightarrow x_{i+1}
\rightarrow y_{i+1},$
 be the result of the algorithm, where $i\geq 0$, the vertex $x_j$
lies on a component of 
$F$
which is
a path of length 
two for each $j\leq i$, the vertex $x_{i+1}$ lies on a component of
$F$ which is a path of 
length at least 4,  and $d_F(y_{i+1})=2$. 
We define a new pseudo path factor $F'$ by setting $V(F')=V(F)$ and \medskip

 \centerline{$E(F')=(E(F)\setminus \{x_jy_j : 
j=1,\dots,i, i+1\})\cup 
\{y_{j-1}x_j : j=1,\dots ,i, i+1\}$. }
\bigskip
\noindent  
Clearly, $V_{F'}=V_F\cup \{y_0\}$ and the proof of Theorem \ref{th:pseudo} is complete.
\end{proof}

Theorems \ref{th:pseudofactor} and \ref{th:pseudo} imply the following
theorem:

\begin{theorem}
\label{th:factor}
Every simple $(3,4)$-biregular bigraph  has a  path factor  such that
the endpoints of each 
path have degree 3.
\end{theorem}

\end{document}